\documentclass[11pt]{article}
\usepackage{amssymb}
\newtheorem{precor}{{\bf Corollary}}

\newenvironment{cor}{\begin{precor}{\hspace{-0.5
               em}{\bf.\ }}}{\end{precor}}
\newtheorem{precon}{{\bf Conjecture}}

\newenvironment{con}{\begin{precon}{\hspace{-0.5
               em}{\bf.\ }}}{\end{precon}}
\newtheorem{predefin}{{\bf Definition}}

\newenvironment{defin}[1]{\begin{predefin}{\hspace{-0.5
                   em}{\bf.\ }}{\rm #1}\hfill{$\spadesuit$}}{\end{predefin}}
\newtheorem{preexm}{{\bf Example}}

\newtheorem{preappl}{{\bf Application}}

\newtheorem{prelem}{{\bf Lemma}}

\newenvironment{lem}{\begin{prelem}{\hspace{-0.5
               em}{\bf.\ }}}{\end{prelem}}
\newtheorem{preproof}{{\bf Proof.\ }}

\newenvironment{proof}[1]{\begin{preproof}{\rm
               #1}\hfill{$\blacksquare$}}{\end{preproof}}
\newtheorem{prethm}{{\bf Theorem}}

\newenvironment{thm}{\begin{prethm}{\hspace{-0.5
               em}{\bf.\ }}}{\end{prethm}}
\newtheorem{prealphthm}{{\bf Theorem}}

\newenvironment{alphthm}{\begin{prealphthm}{\hspace{-0.5
               em}{\bf.\ }}}{\end{prealphthm}}

\newtheorem{prealphlem}{{\bf Lemma}}

\newenvironment{alphlem}{\begin{prealphlem}{\hspace{-0.5
               em}{\bf.\ }}}{\end{prealphlem}}

\newtheorem{prepro}{{\bf Proposition}}

\newtheorem{preprb}{{\bf Problem}}

\newtheorem{preapp}{{\bf Application}}

\newtheorem{prequ}{{\bf Question}}

\newenvironment{qu}{\begin{prequ}{\hspace{-0.5
               em}{\bf.\ }}}{\end{prequ}}
\def\conct[#1,#2]{\mbox {${#1} \leftrightarrow {#2}$}}
\def\dconct[#1,#2]{\mbox {${#1} \rightarrow {#2}$}}
\def\deg[#1,#2]{\mbox {$d_{_{#1}}(#2)$}}
\def\mindeg[#1]{\mbox {$\delta_{_{#1}}$}}
\def\maxdeg[#1]{\mbox {$\Delta_{_{#1}}$}}
\def\outdeg[#1,#2]{\mbox {$d_{_{#1}}^{^+}(#2)$}}
\def\minoutdeg[#1]{\mbox {$\delta_{_{#1}}^{^+}$}}
\def\maxoutdeg[#1]{\mbox {$\Delta_{_{#1}}^{^+}$}}
\def\indeg[#1,#2]{\mbox {$d_{_{#1}}^{^-}(#2)$}}
\def\minindeg[#1]{\mbox {$\delta_{_{#1}}^{^-}$}}
\def\maxindeg[#1]{\mbox {$\Delta_{_{#1}}^{^-}$}}
\def\isdef{\mbox {$\ \stackrel{\rm def}{=} \ $}}
\def\dre[#1,#2,#3]{\mbox {${\cal E}^{^{#3}}(#1,#2)$}}
\def\var[#1,#2]{\mbox {${\rm Var}_{_{#1}}(#2)$}}
\def\ls[#1]{\mbox {$\xi^{^{#1}}$}}
\def\hom[#1,#2]{\mbox {${\rm Hom}({#1},{#2})$}}
\def\onvhom[#1,#2]{\mbox {${\rm Hom^{v}}(#1,#2)$}}
\def\onehom[#1,#2]{\mbox {${\rm Hom^{e}}(#1,#2)$}}
\def\core[#1]{\mbox {$#1^{^{\bullet}}$}}
\def\cay[#1,#2]{\mbox {${\rm Cay}({#1},{#2})$}}
\def\sch[#1,#2,#3]{\mbox {${\rm Sch}({#1},{#2},{#3})$}}
\def\cays[#1,#2]{\mbox {${\rm Cay_{s}}({#1},{#2})$}}
\def\dirc[#1]{\mbox {$\stackrel{\rightarrow}{C}_{_{#1}}$}}
\def\cycl[#1]{\mbox {${\bf Z}_{_{#1}}$}}

\begin{document}
\footnotetext[1]{This paper is partially supported by
Shahid Beheshti University.}
\begin{center}
{\Large \bf Circular Coloring and Mycielski Construction}\\

{\bf Meysam Alishahi and Hossein Hajiabolhassan}\\
{\it Department of Mathematical Sciences}\\
{\it Shahid Beheshti University, G.C.}\\
{\it P.O. Box {\rm 1983963113}, Tehran, Iran}\\
{\tt m\_alishahi@sbu.ac.ir}\\
{\tt hhaji@sbu.ac.ir}\\ \ \\
\end{center}
\begin{abstract}
\noindent In this paper, we investigate circular chromatic number
of Mycielski construction of graphs. It was shown in
\cite{MR2279672} that $t^{{\rm th}}$ Mycielskian of the Kneser
graph $KG(m,n)$ has the same circular chromatic number and
chromatic number provided that $m+t$ is an even integer. We prove
that if $m$ is large enough, then
$\chi(M^t(KG(m,n)))=\chi_c(M^t(KG(m,n)))$ where $M^t$ is $t^{{\rm
th}}$ Mycielskian. Also, we consider the generalized Kneser graph
$KG(m,n,s)$ and show that there exists a threshold $m(n,s,t)$
such that $\chi(M^t(KG(m,n,s)))=\chi_c(M^t(KG(m,n,s)))$ for
$m\geq m(n,s,t)$.\\

\end{abstract}
\section{Introduction}
In this section we elaborate on some basic definitions and facts
that will be used throughout the paper. Throughout this paper we
only consider finite graphs. A {\it homomorphism} $f: G
\longrightarrow H$ from a graph $G$ to a graph $H$ is a map $f:
V(G) \longrightarrow V(H)$ such that $uv \in E(G)$ implies
$f(u)f(v) \in E(H)$. Also, the symbol $\hom[G,H]$ is used to
denote the set of all homomorphisms from $G$ to $H$.

For a given graph $G$, the notation $g(G)$ stands for the girth
of graph $G$. We denote the neighborhood of a vertex $v\in V(G)$
by $N(v)$ and $N[v]$ stands for the closed neighborhood of $v$,
i.e., $N[v]=N(v) \cup \{v\}$. We denote by $[m]$ the set $\{1, 2,
\ldots, m\}$, and denote by ${[m] \choose n}$ the collection of
all $n$-subsets of $[m]$. The {\it Kneser graph} $KG(m,n)$ is the
graph with vertex set ${[m] \choose n}$, in which $A$ is
connected to $B$ if and only if $A \cap B = \emptyset$. It was
conjectured by Kneser \cite{MR00685361} in 1955, and proved by
Lov\'{a}sz \cite{MR514625} in 1978, that $\chi(KG(m,n))=m-2n+2$.
A subset $S$ of $[m]$ is called $2$-{\it stable} if  $2 \le |x-y|
\le m-2$ for all distinct elements $x$ and $y$ of $S$.  The {\it
Schrijver graph} $SG(m,n)$ is the subgraph of $KG(m,n)$ induced
by all $2$-stable $n$-subsets of $[m]$. It was proved by Schrijver
\cite{MR512648} that $\chi(SG(m,n))=\chi(KG(m,n))$ and every
proper subgraph of $SG(m,n)$ has a chromatic number smaller than
that of $SG(m,n)$. The {\it fractional chromatic number},
$\chi_{_{f}}(G)$, of a graph $G$ is defined as
$$\chi_{_{f}}(G) \isdef \inf \{\frac{m}{n} \ \ | \ \ \hom[G,KG(m,n)]
\not = \emptyset \}.$$

The so-called {\it generalized Kneser graphs} are generalized from
the Kneser graphs in a natural way. Let $m\geq 2n$ be positive
integers. The  generalized Kneser graphs $KG(m, n, s)$ is the
graph whose vertex-set is the set of $n$-subsets of $[m]={1,
2,\ldots, m}$ where two $n$-subsets $A$ and $B$ are joined by an
edge if $|A\cap B|\leq s$. Some properties of generalized Kneser
graphs has been studied in several papers, see \cite{MR2427759,
MR1468332, MR797510}.

The Erd${\rm \ddot{o}}$s-Ko-Rado theorem \cite{MR0140419} states
that if ${\cal F}$ is an independent set in the Kneser graph
$KG(m,n)$, then $|{\cal F}|\leq{{m-1}\choose {n-1}}$. If $m > 2n$
and equality holds, then $F$ is trivial, i.e., $|\displaystyle
\cap_{F\in{\cal F}}F|=1$. Here is a generalization of the
Erd${\rm \ddot{o}}$s-Ko-Rado theorem.

\begin{alphthm}{\rm \cite{MR771733}}\label{indgen}
Let $n > s$ be non-negative integers. If $m\geq (s+2)(n-s)$ and
${\cal F}$ is an independent set of the generalized Kneser graph
$KG(m,n,s)$, then $|{\cal F}|\leq{{m-s-1}\choose {n-s-1}}$.
\end{alphthm}

Ahlswede and Khachatrian \cite{MR1429238} determined the maximum
size of independent sets for all $m$, by proving their complete
intersection theorem.

Let $G$ be a graph. A $k$-coloring of $G$ is the partition of the
vertex set into independent sets $V_1, V_2, \ldots, V_k$. Also,
$V_1, V_2, \ldots, V_k$ is termed a $k$-coloring of $G$. If $n$
and $d$ are positive integers with $n \geq 2d$, then the {\it
circular complete graph} $K_{n\over d}$ is the graph with vertex
set $\{0,1, \ldots, n-1\}$ in which $i$ is connected to $j$ if
and only if $d \leq |i-j| \leq n-d$. A graph $G$ is said to be
$(n, d)$-colorable if $G$ admits a homomorphism to $K_{n\over
d}$. The {\it circular chromatic number} (also known as the {\it
star chromatic number} \cite{MR968751}) $\chi_{_{c}}(G)$ of a
graph $G$ is the minimum of those ratios $\frac{n}{d}$ for which
${\rm gcd}(n,d)=1$ and such that $G$ admits a homomorphism to
$K_{n\over d}$. It can be shown that one may only consider
onto-vertex homomorphisms \cite{MR1815614}. It is known
\cite{MR968751,MR1815614} that for any graph $G$, $\chi(G) - 1 <
\chi_c(G) \leq \chi(G)$, and hence $\chi(G) = \lceil \chi_c(G)
\rceil$. So $\chi_c(G)$ is a refinement of $\chi(G)$, and
$\chi(G)$ is an approximation of $\chi_c(G)$.

The circular chromatic number of Kneser graphs has been studied
by Johnson, Holroyd, and Stahl \cite{MR1475894}. They proved that
$\chi_c({\rm KG}(m,n))=\chi({\rm KG}(m,n))$ if $m\leq 2n+2$ or
$n=2$, and conjectured that equality holds for all Kneser graphs.

\begin{con}
\label{jhsconj} {\rm \cite{MR1475894}} For all $m \geq 2n+1$,
$\chi_c({\rm KG}(m,n))=\chi({\rm KG}(m,n))$.
\end{con}

This conjecture has been studied in several papers
\cite{MR2340388, MR1983360,MR2197228,MR2279672}. It was proved in
\cite{MR1983360}, if $m$ is sufficiently large, then the
conjecture holds.  Later, it was shown in
\cite{MR2197228,MR2279672}, if $m$ is an even positive integer,
then the conjecture is true.

\begin{alphthm}{\rm \cite{MR1983360}}\label{haji}
If $m \geq 2n^2(n-1)$, then $\chi_c({\rm KG}(m,n)) =\chi({\rm
KG}(m,n))$.
\end{alphthm}

The following notations will be needed throughout the paper. Let
$G$ be a graph with vertex set $\{v_1,v_2,\ldots,v_n\}$. Recall
that the Mycielskian $M(G)$ of $G$ is the graph defined on
$\{v_1,v_2,\ldots,v_n\}\cup\{{v'}_1,{v'}_2,\ldots,{v'}_n\}\cup\{z\}$
with edge set $E(M(G))=E(G)\cup\{{v'}_iv_j:\ v_iv_j\in
E(G)\}\cup\{z{v'}_i:\ i=1,2,\ldots,n \}$.  The vertex ${v'}_i$ is
called the {\it twin} of the vertex $v_i$ (${v}_i$ is also called
the {\it twin} of the vertex ${v'}_i$); and the vertex $z$ is
called the {\it root} of $M(G)$. If there is no ambiguity we
shall always use $z$ as the root of $M(G)$. For $t \geq 2$, let
$M^t(G)\isdef M(M^{t-1}(G))$. Hereafter, for $t\geq 1$, the roots
of $M^t(G)$ are the set of the roots, their twins, the twins of
the twins, and etc., i.e., the set of the roots of $M^t(G)$ is
equal to the vertex set of $M^{t-1}(z)$ where $z$ is the root of
$M(G)$. Mycielski \cite{MR0069494} used this construction to
increase the chromatic number of a graph while keeping the clique
number fixed: $\chi(M(G))=\chi(G)+1$ and $\omega(M(G))=\omega(G)$.

The problem of calculating the circular chromatic number of the
Mycielskian of graphs has been investigated in the literature
\cite{MR1703312, MR2057687, MR2003515, MR1704152, MR2062859,
MR2279672}. It turns out that the circular chromatic number of
$M(G)$ is not determined by the circular chromatic number of $G$
alone. Rather, it depends on the structure of $G$. Even for graph
$G$ with very simple structure, it is still difficult to
determine $\chi_c(M(G))$. The problem of determining the circular
chromatic number of the iterated Mycielskian of complete graphs
was discussed in \cite{MR1703312}. It was conjectured in
\cite{MR1703312} that if $n\geq t+2\geq3$, then
$\chi_c(M^t(K_n))=\chi(M^t(K_n))=n+t$.

It is known that $\chi_c(M^t(G))=\chi(M^t(G))=n+t$ whenever
$t\geq2$ and $m\geq 2^{t-1}+2t-2$, see \cite{MR2062859}. This
improves a similar result in \cite{MR2003515}. In
\cite{MR2279672}, Simonyi and Tardos,  by using algebraic
topology, showed  if $n+t$ is even, then
$\chi_c(M^t(K_n))=\chi(M^t(K_n))=t+n$.

In what follows we are concerned with some results concerning the
circular chromatic number of the Mycielskian of graphs. In the
first section, we set up notations and terminologies and present
some preliminaries.  In the second section, we introduce the
concept of  free coloring of graphs. Section $2$ establishes the
relationship between circular chromatic number and free chromatic
number. Section $3$ contains a generalization of the concept of
free coloring. In section $4$, we introduce some sufficient
conditions for equality of circular chromatic number and chromatic
number of graphs in terms of $(a,b)$-free chromatic number. In
the last section, we consider the generalized Kneser graph
$KG(m,n,s)$ and show that there exists a threshold $m(n,s,t)$
such that $\chi(M^t(KG(m,n,s)))=\chi_c(M^t(KG(m,n,s)))$ for $m\geq
m(n,s,t)$. Also,  We show that if $m\geq
2n^2(n-1)+\min\{2^{t+1}-2,\ 2^t+3\}n-\min\{0,2n-t-3\}$, then
$\chi(M^t(KG(m,n)))=\chi_c(M^t(KG(m,n)))$ where $M^t$ is $t^{{\rm
th}}$ Mycielskian.
\section{Free chromatic number}

In this section, we introduce the concept of free coloring and
free chromatic number of graphs. We show that when the free
chromatic number of a graph $G$ is large enough, then
$\chi(G)=\chi_c(G)$. In this regard, we introduce some bounds for
free chromatic number of graphs.

\begin{defin}{For a given graph $G$ an independent set is called a {\it free indepndent
set} of $G$ if it can be extended to at least two distinct
maximal independent sets of $G$.}
\end{defin}

It is a simple matter to check that a vertex of graph is
contained in a free independent set if and only if the graph
obtained by deletion of this vertex and its neighborhood is a
non-empty graph. Also, it is easily seen that for a graph $G$ an
independent set $F\subseteq V(G)$ is a free independent set of $G$
if and only if there exists an edge $uv\in E(G)$ such that
$(N(u)\cup N(v))\cap F=\emptyset$.

\begin{defin}{
Let $G$ be a graph. If for any vertex $v$ of $G$, there exists a
free independent set which contains $v$, then $G$ is called a {\it
free graph}.}
\end{defin}

\begin{defin}{
{\it Free chromatic number} of a graph $G$, denoted by $\phi(G)$,
is the minimum positive integer $t$ such that there exists a
partitioning of the vertices of $G$ such as $V(G)= V_1\cup V_2
\cup \cdots \cup V_t$ where $V_i$ is a free independent set for
any $1\leq i \leq t$. If $G$ is not free, then we define the free
chromatic number of $G$ to be $\infty$.}
\end{defin}

The following simple lemma provides a necessary condition for the
existence of graph homomorphism based on the free chromatic
number of graphs.

\begin{lem}
Let $G$ and $H$ be connected free graphs. If there exists an onto
edge homomorphism from $G$ to $H$, then $\phi(G) \leq \phi(H)$.
\end{lem}

An easy computation shows that if $(n,d)=1$ and $d\geq 2$, then
the free chromatic number of circular complete graph $K_{n\over
d}$ is less than twice of its chromatic number.

\begin{lem}\label{freecircular}
Let $G$ be a graph and $\chi_c(G)=\frac nd$ such that $(n,d)=1$.
If $d\geq2$ or equivalently $\chi_c(G)\neq\chi(G)$, then
$\phi(G)\leq \lceil\chi_c(G)(1+\frac1{d-1})\rceil \leq
2\chi(G)-1$.
\end{lem}
\begin{proof}{
Since the circular chromatic number of $G$ is ${n\over d}$, Hence
there exists a homomorphism from graph $G$ to the circular
complete graph $K_{n\over d}$. It was proved that for any $i$,
the edge between the vertices $i$ and $i+d$ (in $K_{n\over d}$)
should be in the range of homomorphism. Consequently, the inverse
image of any $d-1$ consecutive vertices of $K_{n\over d}$ is a
free independent set of $G$. Hence, the free chromatic number of
graph $G$ is less than or equal to $\lceil{n \over
{d-1}}\rceil=\lceil {n\over d}{d\over {d-1}}\rceil=\lceil
\chi_c(G)(1+{1\over{d-1}})\rceil$. Also,
$\lceil\chi_c(G)(1+\frac1{d-1})\rceil\leq 2\chi(G)-1$ which
completes the proof.}
\end{proof}

The aforementioned lemma provides a sufficient condition for the
equality of chromatic number and circular chromatic number of
graphs. Hence, it can be of interest to have bounds for the free
chromatic number of graphs.

Suppose $G$ is a free graph with $n$ vertices. Let
$\bar{\alpha}(G)$ be the size of the largest free independent set
in $G$. It is obviously that $\phi(G)\geq \frac
{n}{\bar{\alpha}(G)}$. Also, for any edge $e=uv$ in $G$, let
$d(e)$ be the number of vertices that lie in the neighborhood of
$u$ or $v$. Define $d(G)=\displaystyle \min \{d(e): e\in E(G)
\}$. It is obviously that $\bar{\alpha}\leq n-d(G)$. So we have
the next lemma.
\begin{lem}
Let $G$ be a graph with $n$ vertices, then $\phi(G)\geq {n\over
{\bar{\alpha}(G)}} \geq \frac{n}{n-d(G)}$.
\end{lem}

Next theorem gives an upper bound for the free chromatic number of
a graph $G$ in terms of chromatic number and $d(G)$. Also, it is
shown when the girth of a graph $G$ is greater than $4$, then the
difference between its free chromatic number and its chromatic
number is at most $4$.

\begin{thm}\label{girth}
Let $G$ be a free graph. Then $\phi(G)\leq \chi(G)+d(G)$. Also,
if $g(G)\geq5$, then $\phi(G)\leq \chi(G)+4.$
\end{thm}
\begin{proof}{Assume that $e=uv$ is an edge in $G$ such that
$d(e)=d(G)$. Color the vertices of the induced graph $G\setminus
(N(u)\cup N(v))$ by using  at most $\chi(G)$ colors. Assign to
each vertex of $N(u)\cup N(v)$ a new color. One can easily check
that this coloring is a free coloring for $G$. Consequently,
$\phi(G)\leq \chi(G)+d(G)$.

Assume that the girth of $G$ is greater than $4$. First, suppose
$G$ is a tree. In view of the assumption, since $G$ is a free
graph, the diameter of $G$ should be greater than $4$. Let $u_1,
u_2, \ldots,u_d$ be the longest path in $G$. Let $V_1,V_2$ be a
$2$-coloring for $G\setminus N[u_{d-1}]$. Set $V_3\isdef
\{u_{d-1}\}$ and $V_4\isdef N(u_{d-1})$. One can check that
$V_1,V_2,V_3,V_4$ is a free coloring for $G$. Now, suppose that
$G$ is not a tree. Let $C$ be an induced cycle of $G$ such that it
has $g\isdef g(G)$ vertices. Assume that
$V(C)=\{u_1,u_2,\ldots,u_g\}$ and $E(C)=\{u_iu_{i+1}| 1\leq i
\leq g\ ({\rm mod}\ g) \}$. Set $V_1\isdef \{u_1\}$, $V_2\isdef
\{u_2\}$, $V_3\isdef N(u_1)$, $V_4\isdef N(u_2)$. Also, suppose
$V_5,V_6,\ldots,V_{\chi(G)+4}$ is a $\chi(G)$-coloring of
$G\setminus (V_1\cup V_2\cup V_3\cup V_4)$. It is easy seen that
$V_1,V_2,\ldots,V_{\chi(G)+4}$ is a free coloring for $G$.}
\end{proof}
Similarly, if $G$ is a free graph and $g(G)\geq 7$, then set
$V_1\isdef N(u_1)$ and $V_2\isdef N(u_2)$. Moreover, assume that
$V_3,V_4,\ldots,V_{\chi(G)+2}$ is a $\chi(G)$-coloring of
$G\setminus (V_1\cup V_2)$. One can check that
$V_1,V_2,\ldots,V_{\chi(G)+2}$ is a free coloring for $G$. Thus,
$\phi(G)\leq \chi(G)+2.$

Since $d(G)\leq \Delta+\delta$ we have the next corollary.
\begin{cor}
Let $G$ be a free graph. Then $\phi(G)\leq \chi(G)+\Delta+\delta$.
\end{cor}

This question naturally arises in mind that how much free
chromatic number can be larger than chromatic number. Consider
$K_n\times K_m$ which is the Categorical product of the complete
graphs $K_m$ and $K_n$. It is easily seen that the chromatic
number of this graph is $\min\{m,n\}$. On the other hand, every
independent set of size at least $2$ can be extended to a unique
maximal independent set. Hence, the free chromatic number of this
graph is $mn$. In special case, when $n=2$ and $m\geq 2$,
$\phi(K_2\times K_m)=2m=\chi(K_2\times K_m)+d(K_2\times K_m)$.
Therefore, the bound in the previous corollary is sharp.

In view of proof of Theorem \ref{girth}, if $G$ is a graph which
contains two adjacent vertices $u$ and $v$ such that the induced
subgraph on $N(u) \cup N(v)$ is a free graph with girth greater
than $4$, then one can similarly show that $\phi(G)\leq
\chi(G)+4$. Hence, if $G$ contains an induced sparse subgraph,
then Lemma \ref{freecircular} is not fruitful. This leads us to
generalize the definition of free coloring.
\section{Generalization of Free Coloring}

In this section, we generalize the concept of free coloring in
order to show that the circular chromatic number of the ${\rm
t^{th}}$ Mycielskian of the generalized Kneser graph $KG(m,n,s)$
is equal to its chromatic number whenever $m$ is sufficiently
large.

\begin{defin}
{For a free independent set $F$, we say $F$ is supported by the
edge $e=uv$ or $e$ supports $F$ if $F\cap (N(u)\cup
N(v))=\varnothing$. The set of all edges which support $F$ is
noticed by $supp(F)$.}
\end{defin}

Here is a generalization of free chromatic number.

\begin{defin}{Let $a \geq 0$ and $b \geq 1$ be integers. The {\it $(a,b)$-free chromatic number} of a graph $G$,
denoted by $\phi^a_b(G)$, is the minimum natural number $t$ (if
there is not such $t$, we define $\phi^a_b(G)=\infty$) such that
\begin{enumerate}
\item There exists a partition of the vertices of $G$ into independent sets
$V_1,V_2,\ldots, V_t$ where all but at most $a$ of $V_i$'s are
free independent sets. For convince, let $V_i$ be a free
independent sets for $i=1,2,\ldots t-a$.
\item There exist the edges $e_1,e_2,\ldots,e_{t-a}$ such that
for any $1\leq i \leq t-a$, $e_i\in supp(V_i)$  and also for every
vertex $v\in V(G)$ the number of edges among
$e_1,e_2,\ldots,e_{t-a}$ which are incident with $v$ is less than
or equal to $b$.
\end{enumerate}\vspace*{-.78cm}}\end{defin}
If the partition $V(G)= V_1\cup V_2 \cup \cdots \cup V_t$ and
edges $e_1,e_2,\ldots,e_{t-a}$ satisfy the conditions 1 and 2 in
the previous definition where $e_i\in supp(V_i)$, then
$V_1,V_2,\ldots,V_t$ together $e_1,e_2,\ldots,e_{t-a}$ be noticed
as an $(a,b)$-free coloring for $G$. It should be noted there is
no obligation that the edges $e_1,e_2,\ldots,e_{t-a}$ are
distinct.

Here are some elementary properties of $(a,b)$-free chromatic
number of graphs. The bellow lemma can be concluded immediately.

\begin{lem}\label{chain} Let $G$ be a graph and $\alpha(G)$ be the independence number of $G$ .
\begin{enumerate}
\item[{\rm a)}] For any integer $b\geq 1$, $\phi^0_b(G)\geq \phi(G)$. Moreover, $\mathop {\lim }\limits_{b \to \infty }\phi_{b}^0(G)=\phi(G)
$.
\item[{\rm b)}] For any integers $a'\geq a \geq 0$ and $b'\geq b \geq 1$, $\phi^{a}_b(G)\geq
\phi_{b'}^{a'}(G)$.

\item[{\rm c)}] For any integers $a \geq 0$ and $b \geq 1$,
$\phi^{a}_b(G)\geq {{|V(G)|-a\alpha(G)}\over {\bar{\alpha}(G)}}$.
\end{enumerate}
\end{lem}

One can deduce the following lemma whose proof is almost
identical to that of Lemma \ref{freecircular} and the proof is
omitted for the sake of brevity.

\begin{lem}\label{phi2}
If there are some integers $a\geq 0$ and $b\geq 2$ such that
$\phi_b^a(G)\geq 2\chi(G)$, then $\chi_c(G)=\chi(G)$.
\end{lem}
\section{Mycielski construction and circular chromatic number}
In this section we find some relationships between the circular
chromatic number and $(a,b)$-free chromatic number of graphs. The
following lemma was observed in \cite{MR2057687}.
\begin{alphlem}{\rm \cite{MR2057687}}\label{lema}
Let $M(G)$ be the Mycielski construction of $G$. Assume that $z$
is the root. Also, for any $v\in V(G)$, let $v'$ be the twin of
$v$. Suppose $\chi_c(M(G))=\frac {n}{d},\ {\rm gcd}(n,d)=1$, and
$d\geq2$. Then there is a homomorphism  $c\in \hom[M(G),K_{\frac
{n}{d}}]$ such that $c(z)=0$, and $c(v)=c(v')$ if $c(v)\notin
[n-d+1,d-1]\ ({\rm mod}\ n).$
\end{alphlem}

\begin{lem}\label{mmm1}
Let $G$ be a graph. Assume that $a\geq 0$ and $b\geq 1$ are
integers. Then $\phi^a_b(M(G))\geq\phi_{2b}^{a+b}(G)$. Moreover,
for any positive integer $t$, $\phi^{0}_{2}(M^t(G))\geq
\phi_{2^{t+1}}^{2^{t+1}-2}(G)$.
\end{lem}
\begin{proof}{Assume that $V(M(G))=\{v_1,v_2,\ldots,v_n,v_1',v_2',\ldots,v_n',z\}$ where
$V(G)=\{v_1,v_2,\ldots, v_n\}$, $v_i'$ is the twin of $v_i$ and
$z$ is the root of $M(G)$. If $\phi^a_b(M(G))=\infty$, then there
is nothing to prove. Suppose that $\phi^a_b(M(G))=t$. Let
$V_1,V_2,\ldots,V_t$ together $e_1,e_2,\ldots e_{t-a}$ be an
$(a,b)$-free coloring for $M(G)$ where $e_i\in supp(V_i)$. Set
$U_i\isdef V_i\cap V(G)$.

Note that at most $b$ edges of $e_i$'s are incident with $z$.
Without loss of generality, assume that $e_i$, for $i=1,2,\ldots,
t-a-b$, is not incident with $z$. If there is an $1\leq i \leq
t-a-b$ such that $e_i=v_sv_k'$, then we define $e_i'=v_sv_k$,
otherwise, set $e_i'=e_i$. It is easy to see that $e_i'$ supports
$U_i$ for any $1\leq i \leq t-a-b$. Also, it is straightforward to
check that for every vertex $v\in V(G)$ the number of edges among
$e'_1,e'_2,\ldots,e'_{t-a-b}$ which are incident with $v$ is less
than or equal to $2b$. Therefore, $U_1,U_2,\ldots, U_t$ together
$e'_1,e'_2,\ldots e'_{t-a-b}$ is an $(a+b,2b)$-free coloring for
$G$. The other assertion follows by induction on $t$.}
\end{proof}

In view of Lemma \ref{phi2} and the preceding lemma, the following
corollary yields.
\begin{cor}\label{mmm11}
Let $G$ be a graph. If $\chi_c(M^t(G))\neq\chi(M^t(G))$, then
$\phi_{2^{t+1}}^{2^{t+1}-2}(G)\leq 2
\chi(M^t(G))-1=2\chi(G)+2t-1$.
\end{cor}

The aforementioned corollary and the next Lemma can be considered
as generalizations of Lemma \ref{freecircular}.

\begin{lem}\label{mmm2}
Let $G$ be a graph and $t$ be a positive integer. If
$\chi_c(M^t(G))\neq\chi(M^t(G))$, then $
\phi_{2^{t+1}}^{2^{t}+3}(G)\leq 2\chi(M^t(G))-1$.
\end{lem}
\begin{proof}{
Suppose $\chi_c(M^t(G))=\frac {n}{d},\ {\rm gcd}(n,d)=1$, and
$d\geq2$. In view of Corollary \ref{mmm11}, if $t=1$, the
assertion holds. Hence, assume that $t\geq 2$. Consider $M^t(G)$
as the Mycielskian of $H\isdef M^{t-1}(G)$ and let $c\in
\hom[M(H),K_{\frac {n}{d}}]$ which satisfies Lemma \ref{lema}.
Assume that $n=k(d-1)+s$ where $0 \leq s \leq d-2$. Set $V_i
\isdef c^{-1}(\{(i-1)(d-1)+1, (i-1)(d-1)+2,\ldots, i(d-1)\})$ for
$i=1,2,\ldots, k$ and $V_{k+1}\isdef
c^{-1}(\{k(d-1)+1,k(d-1)+2,\ldots, n \})$. It is well-known that
for any $i$, the edge between vertices $i$ and $i+d$ (in
$K_{\frac {n}{d}}$) should be in the range of $c$. Hence, for any
$i\in\{1, 2, \ldots, k+1\}$, the set
$E_i=c^{-1}(\{(i-1)(d-1),i(d-1)+1\})$ is not empty. For any
$i\in\{1, 2, \ldots, k+1\}$, choose an $e_i\in E_{i}$. Note that
for each $1\leq j\leq n$, the number of edges among
$e_1,e_2,\ldots,e_{k+1}$ which are incident with some vertices in
$c^{-1}(j)$ is at most 2 (for $d>2$ this number is at most 1). One
can check that $V_1,V_2,\ldots,V_{k+1}$ together
$e_1,e_2,\ldots,e_{k+1}$ is a $(0,b)$-free coloring ($b\leq 2$)
for $M^t(G)$. Obviously,
$2\chi(M^t(G))-1\geq k+1=\lceil\frac{n}{d-1}\rceil$. We consider two cases:\\
Case I) $d>2$\\
In this case, $V_1,V_2,\ldots,V_{k+1}$ together
$e_1,e_2,\ldots,e_{k+1}$ is a $(0,1)$-free coloring and by Lemma
\ref{mmm1} we have
$$k+1\geq\phi_1^0(M^{t}(G))\geq\phi_2^1(M^{t-1}(G))\geq\cdots\geq\phi^{2^{t}-1}_{2^t}(G).$$
On the other hand, we know that $\phi^{2^{t}-1}_{2^t}(G) \geq
\phi^{2^{t}+3}_{2^{t+1}}(G).$

\noindent  Case II) $d=2$\\
Assume that ${\cal R}=\{y_1,y_2,\ldots,
y_{2^{t-1}-1}\}\cup\{{y'}_1,{y'}_2,\ldots,
{y'}_{2^{t-1}-1}\}\cup\{z\}$ be the roots of $M^{t}(G)$ where
$T=\{y_1,y_2,\ldots, y_{2^{t-1}-1}\}$ are the roots of
$H=M^{t-1}(G)$ and $T'=\{{y'}_1,{y'}_2,\ldots, {y'}_{2^{t-1}-1}\}$
are the twins of the vertices of $T$ ($y'_i$ is the twin of
$y_i$).

Set $U^{t-1}_i\isdef V_i\cap V(M^{t-1}(G))$, for $1\leq i \leq
k+1$. For any vertex $v\in V(M^{t}(G))$, let $n(v)$ denote the
number of edges among $e_1,e_2,\ldots,e_{k+1}$ which are incident
with $v$. The number of edges among $e_1,e_2,\ldots,e_{k+1}$ which
are incident with $z$ is $n(z)$. Without loss of generality,
assume that $e_i$, for $i=1,2,\ldots, k+1-n(z)$, is not incident
with $z$. If there is an $1\leq i \leq k+1-n(z)$ such that
$e_i=v_sv'_k$ where $v_s$ and $v_k$ are the vertices of $H$
($v'_k$ is the twin of $v_k$), then we define $e^{t-1}_i\isdef
v_sv_k$, otherwise, set $e^{t-1}_i\isdef e_i$. It is readily seen
that $e^{t-1}_i$ supports $U^{t-1}_i$ for any $1\leq i \leq
k+1-n(z)$. Inductively, after  iterating the aforementioned
procedure $t-1$ times more, we obtain
$U^{0}_1,U^{0}_2,\ldots,U^{0}_{k+1}$ together
$e^{0}_1,e^{0}_2,\ldots,e^{0}_{k+1-p}$ which is a
$(p,2^{t+1})$-free coloring for $G$. It is a simple matter to
check that $p$ is the number of edges among
$e_1,e_2,\ldots,e_{k+1}$ which are incident with at lest one
vertex of ${\cal R}$.

Now we show that $p\leq 2^t+3$. To see this, let $S\subseteq T$
be the set of vertices whose color (in the coloring $c$) is not in
$[n-1,1]\ ({\rm mod}\ n)$. Define $S'=\{y'_j|\ y_j\in S\}\subseteq
T'$. Note that $c$ satisfies Lemma \ref{lema}, hence, for any
$y_j \in S$, $y_j$ and $y'_j$ have the same color. Therefore, for
every vertex $y_j\in S$, $n(y_j)+n(y'_j)\leq2$. Also, it is easy
to check that the number of edges among $e_1,e_2,\ldots,e_{k+1}$
which are incident with vertices in $(R\setminus S)\cup\{z\}$ is
at most 5. Furthermore, for any vertex $v\in V(M^t(G))$, $n(v)\leq
2$; consequently, $p\leq \displaystyle \sum_{y_i\in T}
n(y_i)+\displaystyle  \sum_{y'_i\in T'} n({y'}_i)+n(z)\leq
2|S|+2|T\setminus S|+5 \leq2^t+3$. }\end{proof}
\subsection{Mycielski Construction of Generalized Kneser Graphs}

Here we investigate the circular chromatic number of the Mycielski
construction of generalized Kneser graphs. Although, the exact
value of $\chi(KG(m,n,s))$ is unknown in general, we show that
$\chi_c(M^t(KG(m,n,s)))=\chi(M^t(KG(m,n,s)))$ whenever $m$ is
large enough.

\begin{thm}
For any fixed integers $n> s\geq 0$ and $t\geq 0$, if $m$ is
large enough, then $\chi_c(M^t(KG(m,n,s)))=\chi(M^t(KG(m,n,s)))$.
\end{thm}
\begin{proof}{
First, we show that $\chi(M^t(KG(m,n,s)))\leq {m\choose
s+1}+t=O(m^{s+1})$. To see this, for every vertex $A\in
V(KG(m,n,s))$, choose an arbitrary subset $B\subseteq A$ of size
$s+1$ and define $c(A)\isdef B$. Obviously, $c$ is a proper
coloring for $KG(m,n,s)$. Therefore, $\chi(KG(m,n,s))\leq
{m\choose s+1}$ which implies that $\chi(M^t(KG(m,n,s)))\leq
{m\choose s+1}+t$. Now we show that $\bar{\alpha}(KG(m,n,s))\leq
{2n \choose s+2}{m-s-2\choose n-s-2}=O(m^{n-s-2})$. Let ${\cal F}$
be a free independent set in $KG(m,n,s)$. Since ${\cal F}$ is a
free independent set, there is an edge $AB\in E(KG(m,n,s))$ such
that ${\cal F}\cup \{A\}$ and ${\cal F}\cup \{B\}$ are still
independent. Thus, for any $F\in {\cal F}$, $|F\cap A|\geq s+1$
and $|F\cap B| \geq s+1$; consequently, $|F\cap (A\cup B)| \geq
s+2$. Now by counting, the  desired inequality holds. If
$\chi_c(M^t(KG(m,n,s)))\neq \chi(M^t(KG(m,n,s)))$, then by
Corollary \ref{mmm11} we have
\begin{equation}\label{eq3}
2(\chi(M^t(KG(m,n,s)))-1\geq
\phi_{2^{t+1}}^{2^{t+1}-2}(KG(m,n,s)).
\end{equation}
In view of Theorem \ref{indgen} and Lemma \ref{chain}(c), when
$m\geq (s+2)(n-s)$ we have
$$\phi^{2^{t+1}-2}_{2^{t+1}}(KG(m,n,s))) \geq {{{m}\choose
{n}}- ({2^{t+1}-2}){m-s-1\choose n-s-1} \over {2n \choose
s+2}{m-s-2\choose n-s-2}}=O(m^{s+2}).$$

Since $\chi(M^t(KG(m,n,s)))\leq {m\choose s+1}+t=O(m^{s+1})$,
there exists a threshold $m(n,s,t)$ such that for $m\geq
m(n,s,t)$, (\ref{eq3}) dose not hold, as desired. }\end{proof}

The chromatic number of the generalized Kneser graph $KG(m,n,1)$
has been specified in \cite{MR797510}. The chromatic number of
$KG(m,n,s)$ remains open for $s\geq 2$. Motivated by the
aforementioned theorem, we propose the following question.

\begin{qu}
Let $m, n,$ and $s$ be non-negative integers where $m > n > s$. Is
it true that $\chi_c(KG(m,n,s))=\chi(KG(m,n,s))$?
\end{qu}
It was proved by Hilton and Milner \cite{MR0219428} that if $X$
is an independent set of ${\rm KG}(m, n)$ of size
$${m-1 \choose n-1}-{m-n-1 \choose n-1}+2,$$
then
$$\displaystyle \bigcap_{A \in X} A = \{i\},$$
for some $i \in [m]$. Therefore, any independent set of size
greater than ${m-1 \choose n-1}-{m-n-1 \choose n-1}+1$ can be
extended to a unique maximum independent set. This leads us to
the following lemma.
\begin{lem}\label{hilton}
 Let $m> 2n$ be positive integers. The size of any free independent
 set in the Kneser graph $KG(m,n)$ is less than or equal to
${m-1 \choose n-1}-{m-n-1 \choose n-1}.$ Also, for any integers
$a\geq 0$ and $b\geq 1$, the $(a,b)$-free chromatic number of the
Kneser graph $KG(m,n)$ is at least
$$\phi^a_b(KG(m,n)) \geq {{m \choose n}-a{m-1 \choose n-1}
\over {{m-1 \choose n-1}-{m-n-1 \choose n-1}}}.$$
\end{lem}
\begin{proof}{Let ${\cal F}$ be a free independent set of $KG(m,n)$.
In view of the Hilton and Milner theorem, we should have $|{\cal
F}| \leq {m-1 \choose n-1}-{m-n-1 \choose n-1}+1$. Since ${\cal
F}$ is a free independent set, it can be extended to two distinct
maximal independent sets. If $|{\cal F}| = {m-1 \choose
n-1}-{m-n-1 \choose n-1}+1$, by considering the Hilton and Milner
theorem, there exist two positive integer $i$ and $j$ such that
$\cap_{A \in X} A = \{i, j\}$. Hence, $|{\cal F}|\leq {m-2
\choose n-2}$. On the other hand, it is easy to check that ${m-2
\choose n-2}\leq {m-1 \choose n-1}-{m-n-1 \choose n-1}$ which is
a contradiction.}
\end{proof}

The next theorem was shown in \cite{MR1983360}.
\begin{alphthm}{\rm \cite{MR1983360}}
For any positive integer $n$, if $m$ is sufficiently large, then
we have $\chi_c(SG(m,n))=\chi(SG(m,n))$.
\end{alphthm}

Here is a generalization of the previous theorem.

\begin{thm}\label{SG}
For any integers $n\geq 1$ and $t \geq 0$, there is a threshold
$m(n,t)$ such that
$\chi_c(M^t(SG(m,n)))=\chi(M^t(SG(m,n)))=m-2n+2+t$ whenever $m\geq
m(n,t)$.
\end{thm}
\begin{proof}{
Let $\chi_c(M^t(SG(m,n)))\neq m-2n+2+t$. Set $k\isdef
\min\{2^{t+1}-2,\ 2^{t}+3\}$. By Corollary \ref{mmm11} and Lemma
\ref{mmm2} we have
\begin{equation}\label{eq1}
\phi_{2^{t+1}}^k(SG(m,n))\leq \phi^{0}_{2}(M^t(SG(m,n)))\leq
2(m-2n+2+t)-1.
\end{equation}
On the other hand, the vertex set of $V(SG(m,n))$ has cardinality
${m-n-1 \choose n-1} \frac{m}{n} =O(m^n)$ (each $2$-stable
$n$-subsets of $[m]$ containing $1$ corresponds to an integral
solution of the equation $x_1+x_2+\cdots + x_n=m$ with $x_i \geq
2$. So there are ${m-n-1 \choose n-1}$, $2$-stable $n$-subsets of
$[m]$ containing $1$). In view of Lemma \ref{hilton}, we have
$$\phi_{2^{t+1}}^{k}(SG(m,n)) \geq {{{m-n-1}\choose {n-1}}\frac
mn-k{m-1\choose n-1} \over {{m-1 \choose n-1}-{m-n-1 \choose
n-1}}}=O(m^2).$$ Therefore, there is a threshold $m(n,t)$ such
that if $m\geq m(n,t)$, then
$$
\phi^{0}_{2}(M^t(SG(m,n)))\geq 2\chi(M^t(SG(m,n)))=2(m-2n+2+t).$$
This contradicts (\ref{eq1}).
 }\end{proof}
Here is a generalization of Theorem \ref{haji}.
\begin{thm}
For any integers $n\geq 1$ and $t \geq 0$, if $m\geq
2n^2(n-1)+\min\{2^{t+1}-2,\ 2^t+3\}n-\min\{0,2n-t-3\}$, then
$\chi_c(M^t(KG(m,n)))=\chi(M^t(KG(m,n)))$.
\end{thm}
\begin{proof}{The case $n=1$ was proved in \cite{MR1983360}, hence, assume that
$n\geq2$. For convince, set $k\isdef \min\{2^{t+1}-2,\ 2^{t}+2\}$.
In view of the proof of Theorem \ref{SG}, it is sufficient to show
that the following inequality holds for $m\geq
2n^2(n-1)+kn-\min\{0,2n-t-3\}$
$$
{m-1\choose n-1}-{m-n-1\choose n-1}\leq {{m\choose
n}-k{m-1\choose n-1}\over 2(m-2n+2+t)}.
$$
By double counting we have
$${m-1 \choose n-1}-{m-n-1
\choose n-1}\leq n{{m-2}\choose{n-2}}.$$ It is straightforward to
check that
$$n{{m-2}\choose{n-2}}\leq {{m\choose
n}-k{m-1\choose n-1}\over 2(m-2n+2+t)}$$ for $m\geq
2n^2(n-1)+kn-\min\{0,2n-t-3\}$. This completes the proof.
}\end{proof}

It was shown in \cite{MR1703312}, if $G$ is a graph with chromatic
number $t$, then $\chi_c(M^{t-1}(G))\leq \chi(M^{t-1}(G))-{1\over
2}$. Also, it was conjectured in \cite{MR1703312} that if $n\geq
t+2\geq3$, then $\chi_c(M^t(K_n))=\chi(M^t(K_n))=n+t$.

\begin{qu}
Let $m,n$, and $t$ be non-negative integers where $m > 2n$ and
$0\leq t\leq m-2n$. Is it true that
$\chi_c(M^t(KG(m,n)))=\chi(M^t(KG(m,n)))$?
\end{qu}


\end{document}